\documentclass[11pt]{amsart}

\usepackage[english]{babel}
\usepackage{ucs}
\usepackage[T1]{fontenc}
\usepackage[utf8x]{inputenc}

\usepackage{url}
\usepackage{dsfont}
\usepackage{enumitem}
\usepackage{graphicx}

\usepackage{times}
\usepackage[mathscr]{euscript}

\usepackage[colorlinks=true]{hyperref}

\newcommand{\dd}{\mathrm{d}}
\newcommand{\ind}[1]{\mathbf{1}_{\{#1\}}}
\newcommand{\R}{\mathbb{R}}
\renewcommand{\Pr}{\mathbb{P}}
\newcommand{\Exp}{\mathbb{E}}

\newcommand{\tp}{\tilde{p}}
\newcommand{\tx}{\tilde{x}}
\newcommand{\tX}{\tilde{X}}
\newcommand{\tZ}{\tilde{Z}}
\newcommand{\tbx}{\tilde{\mathbf{x}}}

\newcommand{\uk}{\underline{k}}
\newcommand{\ok}{\overline{k}}

\title{Long time behaviour and mean-field limit of Atlas models}
\author{Julien Reygner}\address{Universit\'e Paris-Est, CERMICS (ENPC), F-77455 Marne-la-Vall\'ee}

\makeatletter
\def\section{\@startsection{section}{1}%
\z@{.7\linespacing\@plus\linespacing}{.5\linespacing}%
{\normalfont\bfseries\centering}}

\def\@settitle{\begin{center}%
  \baselineskip14\p@\relax
    \bfseries
    \LARGE\@title
  \end{center}%
}

\def\@setauthors{%
  \begingroup
  \trivlist
  \centering\footnotesize \@topsep30\p@\relax
  \advance\@topsep by -\baselineskip
  \item\relax
  \andify\authors
  \def\\{\protect\linebreak}%
 {\Large\authors}%
  \endtrivlist
  \endgroup
}

\def\maketitle{\par
  \@topnum\z@ 
  \@setcopyright
  \thispagestyle{firstpage}
  \ifx\@empty\shortauthors \let\shortauthors\shorttitle
  \else \andify\shortauthors
  \fi
  \@maketitle@hook
  \begingroup
  \@maketitle
  \toks@\@xp{\shortauthors}\@temptokena\@xp{\shorttitle}%
  \toks4{\def\\{ \ignorespaces}}
  \edef\@tempa{%
    \@nx\markboth{\the\toks4
      \@nx{\the\toks@}}{\the\@temptokena}}%
  \@tempa
  \endgroup
  \c@footnote\z@
  \def\do##1{\let##1\relax}%
  \do\maketitle \do\@maketitle \do\title \do\@xtitle \do\@title
  \do\author \do\@xauthor \do\address \do\@xaddress
  \do\email \do\@xemail \do\curraddr \do\@xcurraddr
  \do\commby \do\@commby
  \do\dedicatory \do\@dedicatory \do\thanks \do\thankses
  \do\keywords \do\@keywords \do\subjclass \do\@subjclass
}
\makeatother

\begin{document}

\begin{abstract} This article reviews a few basic features of systems of one-dimensional diffusions with rank-based characteristics. Such systems arise in particular in the modelling of financial markets, where they go by the name of Atlas models. We mostly describe their long time and large scale behaviour, and lay a particular emphasis on the case of mean-field interactions. We finally present an application of the reviewed results to the modelling of capital distribution in systems with a large number of agents. \end{abstract}

\maketitle


\section{Introduction} 

The term \emph{Atlas model} was originally introduced in Fernholz' monograph on \emph{Stochastic Portfolio Theory}~\cite{Fer02} to describe a stock market in which the asset prices evolve according to independent and identically distributed processes, except the smallest one which undergoes an additional upward push. The whole growth of the portfolio is thus entirely borne by the stock in lowest position, whence the reference to the Greek Titan holding up the Heavens on his shoulders. By extension, \emph{Atlas models} may broadly refer to any system of one-dimensional quantities whose evolution only depends on their rank in the ordered system. This article aims to shortly review the main features of such systems, with a particular emphasis on their long time and large scale behaviour.


\subsection{Systems of rank-based interacting diffusions}

We shall focus on diffusion processes on the line. In this context, the original Atlas model of size $n \geq 1$ is defined by the system of stochastic differential equations
\begin{equation}\label{eq:atlas}
  \forall i \in \{1, \ldots, n\}, \qquad \dd X^{i,n}_t = n\gamma \ind{X^{i,n}_t = \min_{1 \leq j \leq n} X^{j,n}_t} \dd t + \dd W^{i,n}_t,
\end{equation}
where $\gamma>0$ and the processes $(W^{1,n}_t)_{t \geq 0}, \ldots, (W^{n,n}_t)_{t \geq 0}$ are independent standard Brownian motions. It is a particular case of a system of \emph{rank-based interacting diffusions}, the generic form of which writes
\begin{equation}\label{eq:rb}
  \forall i \in \{1, \ldots, n\}, \qquad \dd X^{i,n}_t = \sum_{k=1}^n \ind{X^{i,n}_t = X^{(k),n}_t} b^{k,n} \dd t + \sum_{k=1}^n \ind{X^{i,n}_t = X^{(k),n}_t} \sigma^{k,n} \dd W^{i,n}_t,
\end{equation}
where $b^{1,n}, \ldots, b^{n,n}$ are given \emph{growth rate} coefficients, the \emph{diffusion} coefficients $\sigma^{1,n}, \ldots, \sigma^{n,n}$ are assumed not to vanish, and for all $t \geq 0$, $X^{(1),n}_t \leq \cdots \leq X^{(n),n}_t$ refer to the order statistics of the vector $(X^{1,n}_t, \ldots, X^{n,n}_t)$. When $X^{1,n}_t, \ldots, X^{n,n}_t$ describe the positions of a system of $n$ particles evolving on the line, the interpretation of~\eqref{eq:rb} is straightforward: the particle with $k$-th position in the ranked system has constant drift $b^{k,n}$ and diffusion coefficient $\sigma^{k,n}$, until it collides with one of its neighbouring particles with which it exchanges its drift and diffusion coefficients. In the context of Stochastic Portfolio Theory, rank-based interacting diffusions are also known as \emph{first-order models}, and we refer to Fernholz' book~\cite{Fer02}, Banner, Fernholz and Karatzas' seminal article~\cite{BanFerKar05}, as well as Fernholz and Karatzas' updated review~\cite{FerKar09} for a detailed introduction to this field. 


\subsection{Relation with reflected Brownian motions}

Two related processes play a central role in the study of~\eqref{eq:rb}:
\begin{enumerate}
  \item the order statistics $(X^{(1),n}_t, \ldots, X^{(n),n}_t)_{t \geq 0}$, which takes its values in the polyhedron $D_n := \{(x^1, \ldots, x^n) \in \R^n : x^1 \leq \cdots \leq x^n\}$ and satisfies 
  \begin{equation}\label{eq:os}
    \forall k \in \{1, \ldots, n\}, \qquad \dd X^{(k),n}_t = b^{k,n} \dd t + \sigma^{k,n} \dd \beta^{k,n}_t + \frac{1}{2} \dd L^{k-1,k}_t - \frac{1}{2} \dd L^{k,k+1}_t,
  \end{equation}
  where the processes $(\beta^{1,n}_t)_{t \geq 0}, \ldots, (\beta^{n,n}_t)_{t \geq 0}$ are independent standard Brownian motions, and $(L^{k-1,k}_t)_{t \geq 0}$ denotes the local time at $0$ of the nonnegative semimartingale $(X^{(k),n}_t-X^{(k-1),n}_t)_{t \geq 0}$ --- we take the obvious convention that $L^{0,1}_t = L^{n,n+1}_t \equiv 0$;
  \item the \emph{gap process} $(Z^{1,n}_t, \ldots, Z^{n-1,n}_t)_{t \geq 0}$, defined by 
  \begin{equation}\label{eq:gap}
    \forall k \in \{1, \ldots, n-1\}, \qquad Z^{k,n}_t := X^{(k+1),n}_t - X^{(k),n}_t,
  \end{equation}
  and which takes its values in the \emph{orthant} $O_n := [0,+\infty)^{n-1}$.
\end{enumerate}
Both these processes can be seen as \emph{reflected} Brownian motions with constant drift vector and diffusion matrix, respectively with normal reflection on $\partial D_n$~\cite{Tan79} and oblique reflection on $\partial O_n$~\cite{Wil87}. This remark allows to establish connections between rank-based interacting diffusions and dynamics arising in the study of spin glasses models called \emph{systems of competing particles}~\cite{RuzAiz05, ArgAiz07, Shk09}, as well as with several models of queuing systems~\cite{HarWil87:stoch, HarWil87:AP, Wil95}.


\subsection{Generalisations}

Possible generalisations of the process defined by~\eqref{eq:rb}, on which we shall not elaborate, include rank-based models driven by L\'evy processes~\cite{Shk11, Sar16:arxiv}, and \emph{second-order} models, also called \emph{hybrid} Atlas models, where the drift depends on both the rank and the index of a particle~\cite{IchFerKarBanPap11, FerIchKar13:AF}. While this article focuses on systems with a finite but possibly large number of particles, countably infinite systems with rank-based evolution were also considered~\cite{PalPit08, ChaPal10, SarTsa16}: the order statistics of such infinite systems can be seen as a generalisation of Harris' Brownian motion~\cite{Har65}, and precise estimates on the fluctuations of the bottom particle were recently obtained~\cite{DemTsa15, HerJarVal15}. Finally, replacing the coefficients $1/2$ and $1/2$ in front of the local times in~\eqref{eq:os} with different weights leads to systems with \emph{asymmetric collisions}~\cite{FerIchKar13:SPA, KarPalShk16, Sar16:AIHP} which are also of interest in nonequilibrium statistical physics~\cite{FerSpoWei15:AAP, FerSpoWei15:EJP}.


\subsection{Outline of the article}

The article is organised as follows. In Section~\ref{s:lt}, we describe the long time behaviour of the solution to~\eqref{eq:rb} for a fixed number of particles and any choice of growth rate and diffusion coefficients. In Section~\ref{s:ls}, we restrict ourselves to coefficients describing mean-field interactions between the particles, and discuss the propagation of chaos phenomenon which occurs when the number of particles grows to infinity. In Section~\ref{s:cd}, an application to the modelling of capital distribution in large systems is detailed. More references regarding the various aspects that we address are given throughout the text.


\section{Long time behaviour for a finite number of particles}\label{s:lt}


\subsection{Weak and strong solutions, multiple collisions}\label{ss:sol} 

In this section, we consider systems of the form of~\eqref{eq:rb} for a fixed number of particles $n \geq 1$. In general, it is not obvious that such systems are well-posed: seen as a stochastic differential equation in $\R^n$, \eqref{eq:rb} has piecewise constant (and therefore in general discontinuous) drift and diffusion coefficients. It was proved by Bass and Pardoux~\cite{BasPar87} that as soon as the diffusion coefficients $\sigma^{1,n}, \ldots, \sigma^{n,n}$ do not vanish, then~\eqref{eq:rb} admits a weak solution, which is global in time and unique in law. To the best of the author's knowledge, strong existence and pathwise uniqueness are not entirely solved yet: it was proved that there exists a unique strong solution up to the first triple collision~\cite{IchKarShk13}, and that such a collision occurs if and only if the sequence $(\sigma^{1,n})^2, \ldots, (\sigma^{n,n})^2$ fails to be concave~\cite{Sar15}; but it does not seem to be known whether the strong solution continues to exist after the collision. We refer to~\cite{IchKar10, IchSar16, BruSar16} for an in-depth study of multiple collisions.


\subsection{The global stability condition} 

An important remark regarding the long time behaviour of the particle system is that the random vector $(X^{1,n}_t, \ldots, X^{n,n}_t) \in \R^n$ defined by~\eqref{eq:rb} cannot converge to an equilibrium distribution, as the centre of mass $\bar{X}^n_t := \frac{1}{n} \sum_{i=1}^n X^{i,n}_t$ is easily seen to be a Brownian motion with drift $\bar{b}^n := \frac{1}{n} \sum_{k=1}^n b^{k,n}$ and therefore does not have a limit when $t$ goes to infinity. In spite of this lack of ergodicity, a stationary behaviour can nonetheless be observed on processes that are not sensitive to the motion of the centre of mass; for instance, the gap process defined in~\eqref{eq:gap} above, or the \emph{centered} system $(\tX^{1,n}_t, \ldots, \tX^{n,n}_t)_{t \geq 0}$ defined by $\tX^{i,n}_t := X^{i,n}_t - \bar{X}^n_t$. The latter is nothing but the orthogonal projection of the original particle system onto the hyperplane
\begin{equation}
  M_n := \{(\tx^1, \ldots, \tx^n) \in \R^n : \tx^1 + \cdots + \tx^n = 0\},
\end{equation}
and we now focus on its ergodic properties, which were mostly described by Pal and Pitman~\cite{PalPit08} and Jourdain and Malrieu~\cite{JouMal08} --- see also~\cite{PalShk14, IchPalShk13} for further concentration of measure estimates and convergence rates.

We first assume that $(\sigma^{1,n})^2 = \cdots = (\sigma^{n,n})^2 = \sigma^2$. Introducing the continuous and piecewise affine function
\begin{equation}
  V^n(x^1, \ldots, x^n) = -\sum_{k=1}^n b^{k,n}x^{(k)}, \qquad x^{(1)} \leq \cdots \leq x^{(n)},
\end{equation}
on $\R^n$, it is readily seen that the drift of $X^n_t = (X^{1,n}_t, \ldots, X^{n,n}_t)$ in~\eqref{eq:rb} coincides with $-\nabla V(X^n_t)$, where the gradient is taken in the distributional sense. As a consequence, the centered process turns out to be ergodic if and only if 
\begin{equation}\label{eq:tZn}
  \tZ^n := \int_{\tbx \in M_n} \exp\left(-\frac{2}{\sigma^2}V(\tbx)\right)\dd\tbx < +\infty, 
\end{equation}
where $\dd\tbx$ refers to the Lebesgue measure on the linear space $M_n$. In this case, the unique stationary probability measure $\tp^n_\infty$ of the process has the density $(\tZ^n)^{-1}\exp\left(-\frac{2}{\sigma^2}V(\tbx)\right)$ with respect to $\dd\tbx$. Notice that the marginal distributions of this measure have exponential tails.

Little algebra shows that~\eqref{eq:tZn} holds if and only if the coefficients $b^{1,n}, \ldots, b^{n,n} \in \R$ satisfy 
\begin{equation}\label{eq:gs}
  \forall m \in \{1, \ldots, n-1\}, \qquad \frac{1}{m} \sum_{k=1}^m b^{k,n} > \frac{1}{n-m} \sum_{k=m+1}^n b^{k,n}.
\end{equation}
The latter condition is called the \emph{global stability condition}, and has the natural interpretation that for any partition of the system into a group of $m$ \emph{leftmost} particles and a group of $n-m$ \emph{rightmost} particles, the average drift of the first group is required to be larger than the average drift of the second group. This ensures the stability of the whole system around its centre of mass. Under the stronger condition that $b^{1,n} > \cdots > b^{n,n}$, it was proved in~\cite{JouMal08} that the law of $(\tX^{1,n}_t, \ldots, \tX^{n,n}_t)$ converges to $\tp^n_\infty$ at an exponential rate. Anticipating on the study of large systems addressed in the next section, let us also mention that this rate is uniform in $n$ under suitable assumptions on the array of coefficents $\{b^{k,n}, 1 \leq k \leq n\}$.

Under the stationary measure $\tp^n_\infty$, the gaps between consecutive particles turn out to be independent and exponentially distributed. That the stationary measure of the gap process has a product-of-exponential form can in fact be directly checked from the theory of reflected Brownian motions~\cite{Wil87,HarWil87:AP} and holds even when the sequence $(\sigma^{1,n})^2, \ldots, (\sigma^{n,n})^2$ is not constant but satisfies the so-called \emph{skew-symmetry condition} that $(\sigma^{k+1,n})^2 - (\sigma^{k,n})^2$ be constant. For completely arbitrary choices of positive coefficients $(\sigma^{1,n})^2, \ldots, (\sigma^{n,n})^2$, the global stability condition~\eqref{eq:gs} remains necessary and sufficient for the existence and uniqueness of a stationary probability measure for the both the centered system and the gap process~\cite{BanFerKar05,IchFerKarBanPap11,JouRey14}. However, these stationary measures are generally no longer explicit.


\subsection{Formation of clouds}

We finally describe the situation in which the global stability condition~\eqref{eq:gs} does not hold. Our discussion is mostly based on the study carried out in~\cite{JouRey14,JouRey16:jhde} of the \emph{Sticky Particle Dynamics} introduced by Brenier and Grenier~\cite{BreGre98}, which can be understood as the small noise limit of~\eqref{eq:rb}. For each $\ell \in \{1, \ldots, n\}$, we call \emph{cluster} of the index $\ell$ the largest subset of consecutive indices $\{\uk, \ldots, \ok\}$ containing $\ell$ and such that either $\uk = \ok$ or
\begin{equation}
  \forall m \in \{\uk, \ldots, \ok-1\}, \qquad \frac{1}{m-\uk+1}\sum_{k=\uk}^m b^{k,n} > \frac{1}{\ok-m}\sum_{k=m+1}^{\ok} b^{k,n}.
\end{equation}
The latter condition is naturally called the \emph{local stability condition}. It is clear that the clusters of two indices are either equal or disjoint, so that one can consider the set of all distinct clusters $C_1, \ldots, C_D$ as a partition of $\{1, \ldots, n\}$ into $D$ consecutive integer intervals. 

In this formalism, the global stability condition~\eqref{eq:gs} is not satisfied if and only if $D \geq 2$. In this case, for $d \in \{1, \ldots, D\}$, the cluster $C_d$ has an average drift defined by 
\begin{equation}
  \bar{b}^{d,n} := \frac{1}{|C_d|} \sum_{k \in C_d} b^{k,n},
\end{equation}
where $|C_d|$ denotes the cardinality of $C_d$. It is easy to see that $\bar{b}^{1,n} \leq \cdots \leq \bar{b}^{D,n}$. We now call \emph{clouds} the unions of clusters having the same average drift, so that the average drifts of consecutive clouds are increasing. Then, coming back to the behaviour of the particle system~\eqref{eq:rb}, it turns out that the clouds drift away from each other, in the sense that after an almost surely finite time, the rightmost particle of a cloud no longer collides with the leftmost particle of the subsequent cloud. Furthermore, within each cloud, the motion of the particles around the centre of mass of the cloud is either ergodic if the cloud is composed of a single cluster, or null recurrent if the cloud contains several clusters. This provides a complete description of the long time behaviour of the dynamics generated by~\eqref{eq:rb}, beyond the globally stable case.


\section{Nonlinear diffusion in the mean-field limit}\label{s:ls}

In this section, we study the asymptotic behaviour of the particle system defined by~\eqref{eq:rb} when the number $n$ of particles is large. Our purpose is to derive macroscopic laws, in the very spirit of hydrodynamic limits from statistical physics. One may expect such macroscopic laws to depend rather heavily on the choice of coefficients $b^{1,n}, \ldots, b^{n,n}$ and $\sigma^{1,n}, \ldots, \sigma^{n,n}$ for different values of $n$. Our standing modelling assumption in this respect is a form of \emph{continuity} of the interaction, inspired by the idea that a particle with rank $k$ should not have a drastically different behaviour from the particles with rank $k-1$ or $k+1$. In the sequel, we therefore assume that there exist continuous functions $b$ and $\sigma$ on $[0,1]$ such that, for all $n \geq 1$, for all $k \in \{1, \ldots, n\}$,
\begin{equation}\label{eq:mf}
  b^{k,n} = b\left(\frac{k}{n}\right), \qquad \sigma^{k,n} = \sigma\left(\frac{k}{n}\right).
\end{equation}
Notice that the original Atlas model~\eqref{eq:atlas} does not fit into this framework, as it would formally correspond to the choice of $b(u) = \gamma \delta_0(u)$, where $\delta_0$ is the Dirac distribution in $0$. Any smooth approximation of this distribution however provides a `smooth' approximation of the Atlas model.

Although some of the results detailed below may hold for non uniformly elliptic diffusion coefficients~\cite{JouRey13,JouRey15}, for the sake of concision we shall always assume that $\sigma^2(u) > 0$ for all $u \in [0,1]$, which we simply denote by $\sigma^2 > 0$. According to the discussion of Subsection~\ref{ss:sol}, this assumption is sufficient to ensure that the particle system~\eqref{eq:rb} is well-defined in the weak sense.


\subsection{Propagation of chaos and nonlinear diffusion process}\label{ss:pc}

With the choice of coefficients~\eqref{eq:mf}, the system of stochastic differential equations~\eqref{eq:rb} rewrites
\begin{equation}\label{eq:rbmf}
  \forall i \in \{1, \ldots, n\}, \qquad \dd X^{i,n}_t = b\left(u_n(t,X^{i,n}_t)\right) \dd t + \sigma\left(u_n(t,X^{i,n}_t)\right) \dd W^{i,n}_t,
\end{equation}
where 
\begin{equation}\label{eq:un}
  u_n(t,x) = \frac{1}{n}\sum_{i=1}^n \ind{X^{i,n}_t \leq x}
\end{equation}
is the empirical distribution function of $X^{1,n}_t, \ldots, X^{n,n}_t$. This expression highlights the fact that the particles only interact through their empirical distribution, which is characteristic of \emph{mean-field systems} in statistical physics. For such systems, a \emph{propagation of chaos phenomenon} is usually observed~\cite{Szn91}: assume that the initial positions $X^{1,n}_0, \ldots, X^{n,n}_0$ are independent and identically distributed according to some probability measure $m$ on the line --- this is the initial \emph{chaos}; then for all $N \geq 1$, it is expected that, in spite of the interaction between the particles induced by the evolution of the system, the collection of processes $(X^{1,n}_t)_{t \geq 0}, \ldots, (X^{N,n}_t)_{t \geq 0}$ behave, in the $n \to +\infty$ limit, as $N$ independent copies of a diffusion process $(X_t)_{t \geq 0}$ --- this is the \emph{propagation of chaos} to positive times. 

An equivalent formulation of the propagation of chaos phenomenon is the fact the empirical measure of the particle system in the space of sample paths, defined by
\begin{equation}\label{eq:nun}
  \nu^n = \frac{1}{n} \sum_{i=1}^n \delta_{(X^{i,n}_t)_{t \geq 0}},
\end{equation}
converges to a deterministic probability measure $P$ on $C([0,+\infty))$. Taking the formal $n \to +\infty$ limit in~\eqref{eq:rbmf} and~\eqref{eq:un}, one can guess that if this convergence holds, then $P$ should be the law of a weak solution $(X_t)_{t \geq 0}$ to the stochastic differential equation
\begin{equation}\label{eq:X}
  \left\{\begin{aligned}
    & \dd X_t = b(u(t,X_t))\dd t + \sigma(u(t,X_t))\dd W_t,\\
    & u(t,x) = \Pr(X_t \leq x).
  \end{aligned}\right.
\end{equation}
The coefficients of this stochastic differential equation depend not only on the value of the unknown $X_t$, but also on its law: this property is called \emph{nonlinearity in McKean's sense}~\cite{McK66,McK67}, owing to the fact that the Fokker-Planck equation describing the evolution of the time marginal distributions $(P_t)_{t \geq 0}$ of $P$ is nonlinear. Indeed, the latter equation writes
\begin{equation}
  \partial_t P_t = \frac{1}{2} \partial_{xx}\left(\sigma^2(u)P_t\right) - \partial_x \left(b(u)P_t\right),
\end{equation}
with $P_t = \Pr(X_t \in \cdot) = \partial_x u(t,\cdot)$ in the distributional sense. It is remarkable that integrating this equation in the space variable yields the closed evolution equation
\begin{equation}\label{eq:cl}
  \partial_t u = \partial_{xx} A(u) - \partial_x B(u), \qquad A(u) := \frac{1}{2}\int_{v=0}^u \sigma^2(v)\dd v, \quad B(u) := \int_{v=0}^u b(v)\dd v,
\end{equation}
for the distribution function $u(t,x)$. The latter equation is a one-dimensional scalar \emph{conservation law}, with general flux function $B$ and nonlinear viscosity function $A$.

Propagation of chaos results for rank-based models with mean-field coefficients were obtained by Bossy and Talay~\cite{BosTal96,BosTal97}, Jourdain~\cite{Jou97,Jou00:SPA,Jou00:MCAP,Jou02:AAP,Jou06}, and Shkolnikov~\cite{Shk12}, with the probabilistic study of conservation laws of the form~\eqref{eq:cl} as a main motivation. Perhaps the most recent result in this line, extracted from~\cite{JouRey13}, states in particular that under the assumptions that $b$ and $\sigma^2$ be continuous with $\sigma^2>0$, and that $m$ have a finite first order moment, then:
\begin{itemize}
  \item there exists a unique weak solution to the nonlinear stochastic differential equation~\eqref{eq:X} with $X_0 \sim m$;
  \item its law $P$ in the space of sample paths is the limit in probability of the empirical measure $\nu^n$ defined by~\eqref{eq:nun};
  \item the distribution function $u(t,x)$ is the unique solution, in an appropriate weak sense, of the conservation law~\eqref{eq:cl} complemented with the initial condition $u(0,x) = m((-\infty,x])$.
\end{itemize}
The convergence of $\nu^n$ to $P$ can be seen as a Law of Large Numbers. A corresponding Central Limit Theorem was proved by Jourdain~\cite{Jou00:MCAP} and very recently generalised by Kolli and Shkolnikov~\cite{KolShk16} who derived a stochastic partial differential equation related with the fluctuations of $\nu^n$ around $P$. Finally, a Large Deviation Principle in the spirit of the Dawson-G\"artner theory was obtained by Dembo, Shkolnikov, Varadhan and Zeitouni~\cite{DemShkVarZei16}.


\subsection{Evolution of the nonlinear process}

We now focus on the description of the evolution of the nonlinear diffusion process $(X_t)_{t \geq 0}$ defined by~\eqref{eq:X}. We first recall the main features of the particle system, described in Section~\ref{s:lt}, adapted to the case of mean-field coefficients~\eqref{eq:mf}:
\begin{itemize}
  \item[(i)] the centre of mass is a Brownian motion with drift $\bar{b}^n = \frac{1}{n}\sum_{k=1}^n b(\frac{k}{n})$ (and variance of order $\frac{1}{n}$);
  \item[(ii)] the centered system converges to an equilibrium measure if and only if the global stability condition $\frac{1}{m} \sum_{k=1}^m b(\frac{k}{n}) > \frac{1}{n-m} \sum_{k=m+1}^n b(\frac{k}{n})$, for all $m \in \{1, \ldots, n-1\}$, is satisfied.
\end{itemize}

Owing to the propagation of chaos, the centre of mass of the particle system is expected to converge to the expectation of the nonlinear diffusion process. The latter satisfies
\begin{equation}
  \forall t \geq 0, \qquad \Exp[X_t] = \Exp[X_0] + \int_{s=0}^t \Exp[b(u(s,X_s))]\dd s.
\end{equation}
Since $u(s,\cdot)$ is the distribution function of $X_s$, the random variable $u(s,X_s)$ is uniformly distributed on $[0,1]$ as soon as $u(s,\cdot)$ is continuous on $\R$, which actually holds at least $\dd s$-almost everywhere under quite general assumptions. As a consequence,
\begin{equation}
  \int_{s=0}^t \Exp[b(u(s,X_s))]\dd s = \bar{b}t, \qquad \bar{b} := \int_{v=0}^1 b(v) \dd v = \lim_{n \to +\infty} \bar{b}^n,
\end{equation}
so that the `centre of mass' of the nonlinear diffusion process travels on the line at constant speed $\bar{b}$.

Let us denote by $\tX_t = X_t - \bar{b}t$ the fluctuation of $X_t$ around $\bar{b}t$. By construction, this process has a constant expectation. Assume that one can find a stationary probability distribution for it, and denote by $\phi$ its distribution function. Then the distribution function of $X_t = \bar{b}t + \tX_t$ is $u(t,x) = \phi(x-\bar{b}t)$, and by the results of Subsection~\ref{ss:pc}, it is a solution to the conservation law~\eqref{eq:cl} which has the shape of a \emph{travelling wave}. Injecting this specific form into~\eqref{eq:cl}, one deduces that $\phi$ must satisfy
\begin{equation}
  \frac{\sigma^2(\phi)}{2}\phi' = B(\phi)-\bar{b}\phi.
\end{equation}
Under the assumption that $\sigma^2>0$, it is quickly seen that $\phi$ exists if and only if the flux function $B$ satisfies the so-called \emph{Oleinik entropy condition}~\cite{Ole59} that $B(u)>\bar{b}u$ for all $u \in (0,1)$. This condition rewrites
\begin{equation}\label{eq:ole}
  \forall u \in (0,1), \qquad \frac{1}{u} \int_{v=0}^u b(v)\dd v > \frac{1}{1-u}\int_{v=u}^1 b(v)\dd v,
\end{equation}
which is the exact continuous equivalent of the global stability condition for the particle system. Under this condition, let us introduce the function
\begin{equation}\label{eq:Psi}
  \Psi(u) = \int_{v=1/2}^u \frac{\sigma^2(v)}{2(B(v)-\bar{b}v)}\dd v
\end{equation}
on $(0,1)$. Then $\Psi$ is a bijection from $(0,1)$ to $\R$, and given a distribution function $\phi$ on the line, $\phi(x-\bar{b}t)$ is a travelling wave solution to~\eqref{eq:cl} if and only if there exists $c \in \R$ such that $\phi(x) = \Psi^{-1}(x+c)$. In other words, all the stationary measures for $(\tX_t)_{t \geq 0}$ are translations of each other, with an explicit distribution function, and to select one of them, one can for instance prescribe the value of its expectation.


\subsection{Stability of travelling waves}\label{ss:stab}

In the previous paragraph, we established that under the Oleinik entropy condition~\eqref{eq:ole}, if one takes $X_0$ distributed according to a probability measure $m$ with distribution function $\phi(x) = \Psi^{-1}(x+c)$ for some $c \in \R$, then $X_t$ writes $\bar{b}t + \tX_t$ where $\tX_t$ describe stationary fluctuations that remain distributed according to $m$. When $m$ is \emph{not} of this form, it is natural to wonder whether $u(t,x)=\Pr(X_t \leq x)$ will approach a function of the form $\phi(x-\bar{b}t)$ when time goes to infinity: from a probability point of view, it is a question of convergence to equilibrium for the (nonlinear) diffusion process $(\tX_t)_{t \geq 0}$; from an analysis point of view, it is a question of stability of travelling waves under perturbations. 

The latter question has been investigated thoroughly since the 50's~\cite{Hop50,IliOle58} due to the physical importance of travelling waves, in particular in the study of hyperbolic conservation laws. In the case of a linear viscosity function $A$ (that is to say a constant function $\sigma^2>0$), it was proved by Osher and Ralston~\cite{OshRal82} and Serre~\cite{Ser02} that as soon as $B$ satisfies the Oleinik entropy condition~\eqref{eq:ole} and $u_0$ is a distribution function on the line such that
\begin{equation}\label{eq:u0phi}
  u_0 - \phi \in L^1(\R), \qquad \int_{x \in \R} (u_0(x) - \phi(x))\dd x = 0,
\end{equation}
then the solution $u$ to the Cauchy problem
\begin{equation}
  \left\{\begin{aligned}
    & \partial_t u = \frac{\sigma^2}{2} \partial_{xx} u - \partial_x B(u),\\
    & u(0,x) = u_0(x),
  \end{aligned}\right.
\end{equation} 
satisfies
\begin{equation}
  \lim_{t \to +\infty} \|u(t,\cdot) - \phi(\cdot-\bar{b}t)\|_{L^1(\R)} = 0.
\end{equation}
This result was generalised to initial data that are $L^1$ but not necessarily bounded perturbations of $\phi$ by Freist\"uhler and Serre~\cite{FreSer98}, and to nonlinear viscosity functions by Gasnikov~\cite{Gas09}. 

In probabilistic terms, the stability result expresses the fact that the fluctuation $\tX_t$, the expectation of which is constant, converges in distribution, when $t \to +\infty$, to the stationary measure with the same expectation --- this is the meaning of the condition~\eqref{eq:u0phi}. This convergence is measured in the $L^1$ distance between distribution functions, which is known to coincide with the Wasserstein distance of order $1$~\cite{BobLed14}. Convergence results in Wasserstein distances of higher order, including the case of a nonlinear viscosity, were obtained in~\cite{JouRey13}. 


\section{Fluid description of capital distribution}\label{s:cd}

This last section is dedicated to the application of the results reviewed in the two previous sections to the modelling of capital distribution in an equity market where the logarithms of the capitalisations of each stock are assumed to evolve according to rank-based interactions of the form~\eqref{eq:rb}.


\subsection{Capital distribution curves}

\begin{figure}
  \centering\includegraphics[width=.6\textwidth]{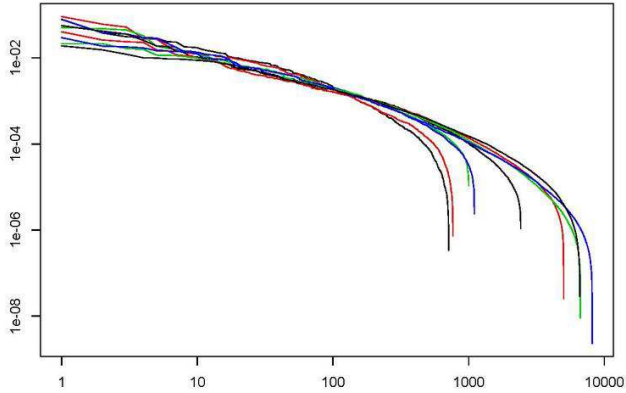}
  \caption{Capital distribution curves for the U.S. stock market. Each curve corresponds to one decade from 1920-1929 to 2000-2009. The size of the market has continuously increased over time, whence the different numbers of stocks for each curve.}
  \label{fig:fer}
\end{figure}

We consider a market of $n$ companies, whose capitalisations at time $t \geq 0$ are denoted by $Y^{1,n}_t, \ldots, Y^{n,n}_t > 0$. Roughly speaking, the capitalisation of a company is the number of its shares times the price of a share: it has to be understood as the total financial value of the company. The \emph{market weight} of the $i$-th company is then defined by
\begin{equation}
  \mu^{i,n}_t = \frac{Y^{i,n}_t}{Y^{1,n}_t + \cdots + Y^{n,n}_t} \in (0,1),
\end{equation}
and simply denotes the proportion of the total wealth held by the company. Writing $\mu^{[1],n}_t \geq \cdots \geq \mu^{[n],n}_t$ for the reverse order statistics of $\mu^{1,n}_t, \ldots, \mu^{n,n}_t$, the log-log plot of the curve $p \mapsto \mu^{[p],n}_t$ is called the \emph{capital distribution curve} at time $t$. By construction, it is nonincreasing. Capital distribution curves for the U.S. stock market are plotted on Figure~\ref{fig:fer}\footnote{This picture, authored by Robert Fernholz, is made available on Wikimedia Commons under the Creative Commons CC0 1.0 Universal Public Domain Dedication.}, which conveys two striking remarks: first, the shape of the curves is extremely stable over time; second, the curves are linear on their leftmost part, which indicates a power law behaviour for the largest market weights.

\subsection{The mean-field capital density}

A cornerstone remark of Fernholz' Stochastic Portfolio Theory~\cite{Fer02} is that, in a market which satisfies certain stability properties, the log-capitalisations $X^{1,n}_t, \ldots, X^{n,n}_t$ defined by
\begin{equation}
  \forall i \in \{1, \ldots, n\}, \qquad X^{i,n}_t := \log Y^{i,n}_t
\end{equation}
are well approximated by the solution of a system of rank-based interacting diffusion of the form~\eqref{eq:rb}. This is the reason why such models are called \emph{first-order} in this context. When one makes the further assumption that the growth rate and diffusion coefficients have the mean-field form~\eqref{eq:mf}, propagation of chaos results become available in order to describe the mean-field limit of such observable quantities as the capital distribution curve. 

In this purpose, we first define the \emph{capital measure} as the probability measure
\begin{equation}
  \pi^n_t := \sum_{p=1}^n \mu^{[p],n}_t \delta_{p/n}
\end{equation}
on $[0,1]$. The interpretation of this measure is straightforward: for all $\alpha \in [0,100]$, $\pi^n_t([0,\frac{\alpha}{100}])$ denotes the proportion of the total capital which is held by the largest $\alpha\%$ companies. The capital distribution curve can then be seen as the histogram (in logarithmic coordinates) of the capital measure, where the ranks of the stocks are rescaled from $\{1, \ldots, n\}$ to the interval $[0,1]$, in order to allow a consistent $n \to +\infty$ limit.

We now recall that the empirical distribution function of the log-capitalisations $X^{1,n}_t, \ldots, X^{n,n}_t$ is denoted by $u_n(t,x)$, $x \in \R$, and write $u_n(t,v)^{-1}$ the empirical quantile of order $v \in [0,1]$, so that for all $p \in \{1, \ldots, n\}$, $X^{[p],n}_t = u_n(t,1-\frac{p-1}{n})^{-1}$. As a consequence of the propagation of chaos result described above, the random probability measure
\begin{equation}
  \pi^n_t = \frac{\sum_{p=1}^n \exp\left(X^{[p],n}_t\right)\delta_{p/n}}{\sum_{q=1}^n \exp\left(X^{[q],n}_t\right)} = \frac{\sum_{p=1}^n\exp\left(u_n(t,1-\frac{p-1}{n})^{-1}\right)\delta_{p/n}}{\sum_{q=1}^n \exp\left(u_n(t,1-\frac{q-1}{n})^{-1}\right)}
\end{equation}
can be proved to converge, when $n \to +\infty$, to the deterministic probability measure with density
\begin{equation}\label{eq:pit}
  \pi_t(v) = \frac{\exp\left(u(t,1-v)^{-1}\right)}{\int_{w=0}^1 \exp\left(u(t,1-w)^{-1}\right)\dd w}
\end{equation}
on $[0,1]$, where $u(t,\cdot)^{-1}$ refers to the pseudo-inverse of the distribution function $u(t,\cdot)$ of the nonlinear diffusion process $X_t$~\cite{JouRey15}. Notice however that exponential integrability conditions on the initial measure $m$ are required for the integral
\begin{equation}
  \int_{w=0}^1 \exp\left(u(t,1-w)^{-1}\right)\dd w = \Exp[\exp(X_t)]
\end{equation}
to be finite. 

Let us call $\pi_t$ the \emph{mean-field capital density}: it gives the $n \to +\infty$ limit of the market weight of groups of companies of size proportional to $n$. However, since the limit of $\pi^n_t$ does not have any atom, the market weight of a single company, or a group of finitely many companies, has vanished. Therefore, we shall say that the mean-field capital density provides a \emph{fluid} description of capital distribution.


\subsection{Stationary behaviour: the Chatterjee-Pal phase transition}

Using the fact that $u(t,\cdot)$ solves the conservation law~\eqref{eq:cl}, a closed-form evolution equation on the pseudo-inverse $u(t,\cdot)^{-1}$ can be derived~\cite{JouRey13}, which then leads to a dynamical fluid description of capital distribution through a closed-form evolution equation on $\pi_t$. Motivated by the remark made above that the shape of capital distribution curves does not seem to vary over time, we only consider the stationary states of this evolution equation, namely the capital density given by the long time limit of $\pi_t$. 

Under the assumptions of Section~\ref{s:ls} regarding the existence and stability of travelling waves for the conservation law~\eqref{eq:cl}, the $t \to +\infty$ limit of~\eqref{eq:pit} is formally obtained by replacing $u(t,\cdot)$ with a travelling wave $\phi(\cdot-\bar{b}t)$, with $\phi(x) = \Psi^{-1}(x+c)$ for some $c \in \R$. One gets $u(t,1-v)^{-1} = \Psi(1-v)+\bar{b}t-c$, so that as soon as $\exp(\Psi(1-v))$ is integrable on $[0,1]$, the terms $\bar{b}t-c$ cancel and $\pi_t$ converges to the stationary capital density given by
\begin{equation}
  \pi_{\mathrm{st}}(v) = \frac{\exp\left(\Psi(1-v)\right)}{\int_{w=0}^1 \exp\left(\Psi(1-w)\right)\dd w}.
\end{equation}
On the other hand, if $\exp(\Psi(1-v))$ is not integrable on $[0,1]$, it can be proved that $\pi_t$ converges to the Dirac mass at $0$. Based on the explicit expression~\eqref{eq:Psi} of $\Psi$, it is observed that whether $\exp(\Psi(1-v))$ is integrable on $[0,1]$ or not directly depends on the coefficients $b$ and $\sigma^2$ as follows:
\begin{itemize}
  \item if $\bar{b}-b(1) > \frac{1}{2}\sigma^2(1)$, then $\exp(\Psi(1-v))$ is integrable on $[0,1]$;
  \item if $\bar{b}-b(1) < \frac{1}{2}\sigma^2(1)$, then $\exp(\Psi(1-v))$ is not integrable on $[0,1]$.
\end{itemize}
In the critical case $\bar{b}-b(1) = \frac{1}{2}\sigma^2(1)$, a more detailed study of $\Psi$ is necessary to determine whether $\exp(\Psi(1-v))$ is integrable on $[0,1]$ or not. 

This `phase transition' was already noted by Chatterjee and Pal~\cite{ChaPal10}, who studied the limit, when $n \to +\infty$, of the stationary distribution of the sequence $(\mu^{[1],n}, \ldots, \mu^{[n],n}, 0, \ldots)$ in the set $\{(\mu^p)_{p \geq 1} : \mu^1 \geq \mu^2 \geq \cdots \geq 0, \sum_{p \geq 1} \mu^p = 1\}$. Since it keeps track of the individual market weights of each company, this approach differs from our \emph{fluid} description. To conclude this article, we discuss the interpretation of the phase transition with both approaches.

\begin{figure}
  \centering\includegraphics[width=.6\textwidth]{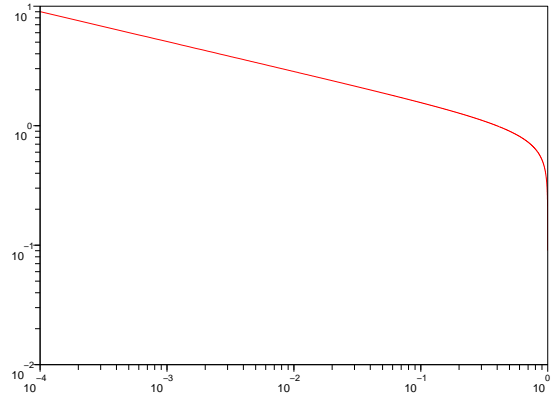}
  \caption{Log-log plot of the stationary density $\pi_{\mathrm{st}}$.}
  \label{fig:capital}
\end{figure}

In the case $\bar{b}-b(1) > \frac{1}{2}\sigma^2(1)$, the stationary capital distribution is described by the density $\pi_{\mathrm{st}}$, meaning that the capital is well spread between all the companies. An instance of a log-log plot of this density is reproduced on Figure~\ref{fig:capital}, which has to be compared with the empirical curves of Figure~\ref{fig:fer}: the shapes are globally similar, and in particular, the power law distribution of largest stocks is captured by our model. Little algebra shows that the slope of the linear part of the curve is equal to $-\sigma^2(1)/2(\bar{b}-b(1))$. In contrast, the Chatterjee-Pal approach only allows to see that the sequence $(\mu^{[1],n}, \ldots, \mu^{[n],n}, 0, \ldots)$ converges to $(0, 0, \ldots)$, which is consistent with the fact that no individual company keeps a `macroscopic' fraction of the total capital.

On the contrary, in the case $\bar{b}-b(1) > \frac{1}{2}\sigma^2(1)$, the stationary capital distribution is the Dirac mass at $0$. This means that all the capital is aggregated by a `microscopic' number of companies, amongst which our fluid description fails to provide a more detailed information on the capital distribution. Yet the  Chatterjee-Pal approach is exactly designed for the study of this case, and it permits to compute the limit distribution of the sequence $(\mu^{[1],n}, \ldots, \mu^{[n],n}, 0, \ldots)$, which turns out to be a Poisson-Dirichlet law~\cite{ChaPal10}. 


\end{document}